\documentclass{amsart}
\usepackage{amssymb,euscript,amsmath, mathrsfs}
\usepackage[dvips]{graphicx}
\usepackage[dvips]{color}

\newcounter{ENUM}

\input xy
\xyoption{all}
\CompileMatrices

\def\<{\langle}
\def\>{\rangle}
\def\0{{{\bf 0}}}

\def\CK{{\mathcal K}}

\def\LL{{\mathbb L}}
\def\tLL{{\tilde {\mathbb L}}}

\def\QQ{{\mathbb Q}}

\def\tR{{\tilde R}}

\def\nubar{{\overline{\nu}}}
\def\rhobar{{\overline{\rho}}}

\newcommand{\Ind}{\operatorname{Ind}}
\newcommand{\cInd}{\operatorname{c-Ind}}

\newcommand{\Rep}{\operatorname{Rep}}

\newcommand{\inv}{\operatorname{inv}}

\def\Hom{\operatorname{Hom}}

\def\End{\operatorname{End}}

\def\GL{\operatorname{GL}}

\def\tr{\operatorname{tr}}
\def\Spec{\operatorname{Spec}}

\newcommand{\margh}[1]{}

\newtheorem{thm}{Theorem}[section]
\newtheorem{prop}[thm]{Proposition}
\newtheorem{lemma}[thm]{Lemma}
\newtheorem{cor}[thm]{Corollary}
\newtheorem{conj}[thm]{Conjecture}

\theoremstyle{definition}
\newtheorem{defn}[thm]{Definition}

\numberwithin{equation}{section}

\begin{document}
\title[Converse theorems]{Converse theorems and the local Langlands correspondence in families}
\author{David Helm}
\author{Gilbert Moss}

\maketitle

\begin{abstract}
We prove a descent criterion for certain families of smooth representations of $\GL_n(F)$ ($F$ a $p$-adic field) in terms of the
$\gamma$-factors of pairs constructed in~\cite{pairs}.  We then use this descent criterion, together with a theory of $\gamma$-factors
for families of representations of the Weil group $W_F$~\cite{galoisgamma}, to prove a series of conjectures, due to the first author,
that give a complete description of the center of the category of smooth $W(k)[\GL_n(F)]$-modules (the so-called ``integral Bernstein center'')
in terms of Galois theory and the local Langlands correspondence.  An immediate consequence is the conjectural ``local Langlands correspondence
in families'' of~\cite{EH-families}.
\end{abstract}

\section{Introduction}

The {\em center} of an abelian category ${\mathcal A}$ is the ring of natural transformations from the identity functor of ${\mathcal A}$
to itself.  It is a commutative ring that acts naturally on every object of ${\mathcal A}$, compatibly with all morphisms of ${\mathcal A}$.

In~\cite{BD}, Bernstein and Deligne study the center of the category of smooth complex representations of a $p$-adic reductive group $G$.  In
particular, they show that such a category is an infinite direct product of full subcategories (called {\em blocks}).  For each such block they
give a concrete and explicit description of the center of each block, and an explicit description of the action of this center on each irreducible
object of the block.

In the context of modular representation theory, however, much less is known.  Paskunas~\cite{paskunas-bernstein} gives a complete description of the center of the
category of finite length representations of $\GL_2(\QQ_p)$ over $p$-adic integer rings, and uses this description to describe the image of the
Colmez functor.  Beyond this, however, the only results currently in the literature are in~\cite{bernstein1}, which describes 
the center of the category of smooth $W(k)[\GL_n(F)]$-modules,
where $F$ is a $p$-adic field and $k$ is an algebraically closed field of characteristic $\ell$ different from $p$.  We will refer to the center
of this category as the {\em integral Bernstein center} in what follows.

Ideally, one would like to have an analogue of the Bernstein-Deligne result for the integral Bernstein center; that is, an explicit description
of the algebra itself, together with its action on simple smooth $W(k)[\GL_n(F)]$-modules.  Results of the first author in~\cite{bernstein1}, show
that although the integral Bernstein center has several nice properties, and can be explicitly described in small examples, it quickly becomes too
complicated (particularly when $n$ is large compared to $\ell$) to admit a description along the lines of Bernstein-Deligne.

Rather than seeking for descriptions of the integral Bernstein center in terms of explicit algebras, it turns out to be better to try to understand
the integral Bernstein center via Galois theory and the local Langlands correspondence.  A first step in this direction was taken in~\cite{bernstein3},
which gives a conjectural description of completions of the integral Bernstein center at maximal ideals in terms of universal framed deformation
rings of Galois representations.  (For the precise statements, see Conjectures 7.5 and 7.6 of~\cite{bernstein3}.)

This conjecture, if established, would have considerable arithmetic implications.  For one thing, the main result of~\cite{bernstein3}
(\cite{bernstein3}, Theorem 7.9)
shows that this conjecture implies the conjectural ``local Langlands correspondence in families'' of~\cite{EH-families}.  Moreover, this conjecture
allows one to construct elements of the integral Bernstein center by considering natural invariants of Galois representations.  For instance, given
an element $w$ of $W_F$, there would be a unique element of the integral Bernstein center that acts on every irreducible admissible representation $\Pi$
in characteristic zero by the scalar $\tr \rho(w)$, where $\rho$ is the representation of $W_F$ attached to $\pi$ via local Langlands.  (Such elements
were first constructed by Chenevier over fields of characteristic zero (\cite{chenevier}, Proposition 3.11), and have also been considered by 
Scholze, again in characteristic zero, in~\cite{scholze-langlands}.)  Ila Varma has noted
that the existence of such elements in the integral Bernstein center is useful for formulating local-global compatibility statements, for instance
for torsion classes in the cohomology of Shimura varieties.

In~\cite{curtis}, the first author further refined this conjecture to a statement about finite type algebras.  In particular, for each block
$\nu$ of the category of smooth representations of $\GL_n$, let $Z_{\nu}$ denote its center.  The construction of section 9 of~\cite{curtis}
gives a pair $(R_{\nu},\rho_{\nu})$, where $R_{\nu}$ is a finite type $W(k)$-algebra, and $\rho_{\nu}: W_F \rightarrow \GL_n(R_{\nu})$ is
an $n$-dimensional representation of $W_F$.  The $k$-points of $\Spec R_{\nu}$ are in bijection with isomorphism classes of
(suitably rigidified) $n$-dimensional representations of $W_F$ over $k$ that correspond, via mod $\ell$ local Langlands, to objects of the block
$\nu$.  Moreover the completion of $R_{\nu}$ at any one of these $k$-points is a universal deformation of the corresponding representation.  Thus
$R_{\nu}$ is a finite type object that interpolates the deformation rings appearing in the conjecture of~\cite{bernstein3}.

The paper~\cite{curtis} then formulates two conjectures relating $Z_{\nu}$ and $R_{\nu}$.  The first of these, which we will henceforth call the ``weak conjecture''
is that there is a map from $Z_{\nu}$ to $R_{\nu}$ satisfying a certain compatibility with the characteristic zero local Langlands correspondence.  The
second conjecture (the so-called ``strong conjecture'') explicitly describes the image of this map.  (We refer the reader to section~\ref{sec:langlands} for a more
detailed discussion of this; in particular the precise statements are given as Conjectures~\ref{conj:weak} and~\ref{conj:strong}, below.)  The strong conjecture
is strictly stronger than the weak conjecture, and implies the earlier conjectures of~\cite{bernstein3} (and thus the ``local Langlands correspondence in families''
of~\cite{EH-families}.)

The main objective of this paper is to establish the ``strong conjecture'' of~\cite{curtis}, and therefore all of the other conjectures.  To do so we
rely heavily on the main result of~\cite{curtis}, which shows that if the ``weak conjecture'' holds for all $\GL_m(F)$ with $m \leq n$, then the
``strong conjecture'' holds for $\GL_n(F)$.  Our approach will be to assume the strong conjecture for $\GL_{n-1}(F)$, and show that it implies
the ``weak conjecture'' for $\GL_n(F)$ (Theorem~\ref{thm:induct}, below.)  
Since both conjectures follow easily from local class field theory for $\GL_1$, a straightforward induction then
gives the desired result.

Our proof of Theorem~\ref{thm:induct} is modeled on Henniart's approach to the $n$ by $n-1$ converse theorem~\cite{converse}.  In particular we require notions
such as Whittaker functions, zeta integrals, and $\gamma$-factors for representations over base rings that can be arbitrary $W(k)$-algebras.  The
correct context for such a theory is the setting of so-called ``co-Whittaker'' modules, developed by the first author in~\cite{bernstein3}.  We recall the
definition and basic theory of these modules in section~\ref{sec:co-whittaker} below.  Of particular importance to us will be the ``universal co-Whittaker modules'';
these live over direct factors $Z_{\nu}$ of the integral Bernstein center, and every co-Whittaker module arises, up to a certain natural equivalence, by base
change from a universal one.  Section~\ref{sec:whittaker} develops some basic ``descent'' criteria for a co-Whittaker module over an algebra $A$; the key
point is that such a module arises by base change from a subalgebra $A'$ if it admits ``sufficiently many'' Whittaker functions that take values in $A'$.
(Theorem~\ref{thm:descent} below.)  Our ultimate goal is to apply this theorem to a certain subalgebra $Z'_{\nu}$ of $Z_{\nu}$; the point is that if we can
find sufficiently many Whittaker functions of the universal co-Whittaker module over $Z_{\nu}$ with values in $Z'_{\nu}$, then the universality will tell us
that the identity map on $Z_{\nu}$ factors through $Z'_{\nu}$, so that $Z'_{\nu}$ is in fact all of $Z_{\nu}$.

The subalgebra $Z'_{\nu}$ in question will be constructed using $\gamma$-factors of pairs.  In~\cite{pairs},
the second author developed a theory of zeta integrals and $\gamma$-factors for co-Whittaker modules that is compatible with the classical theory over
algebraically closed fields of characteristic zero and that satisfies a suitable local functional equation.  We recall this theory in section~\ref{sec:gamma}
below.

In section~\ref{sec:descent} we prove the key technical result of the paper, Corollary~\ref{cor:main descent},
which reformulates our descent theorem for co-Whittaker modules in terms of $\gamma$-factors of pairs.  In particular,
we suppose we have a co-Whittaker $A[\GL_n(F)]$-module $V$, and a subalgebra $A'$ of $A$ such that $A$ is finitely generated as an $A'$ module.
For each direct factor $Z_{\nu'}$ of the integral Bernstein center for $\GL_{n-1}(F)$, we have a corresponding universal co-Whittaker module $W_{\nu'}$.
We show that if the $\gamma$-factors attached to the pairs $V \times W_{\nu'}$ (a priori formal Laurent series with coefficients in $A \otimes Z_{\nu'}$) have
coefficients in $A' \otimes Z_{\nu'}$, then $V$ arises via base change from a co-Whittaker module over $A'$.  

Applying this with $A = Z_{\nu}$ and $V$ the universal co-Whittaker module over $A$ lets us prove that, after completing at a maximal ideal, $Z_{\nu}$
is generated by a suitable set of elements derived from
coefficients of ``universal $\gamma$-factors of pairs''.  Although we do not make direct use of this, we include this in
section~\ref{sec:generators} as a result of independent interest.

In section~\ref{sec:langlands} we prove our main result, by using 
Corollary~\ref{cor:main descent} to show that the strong conjecture for $\GL_{n-1}(F)$ implies the weak conjecture
of $\GL_n(F)$.  Two additional ingredients are necessary.  The first is that it is easy to show that the weak conjecture holds after inverting $\ell$,
using the Bernstein-Deligne description of the Bernstein center for fields of characteristic zero (this is Theorem 10.4 of~\cite{curtis}).  This gives us a map
$Z_{\nu} \rightarrow R_{\nu}[\frac{1}{\ell}]$, compatible with local Langlands; the problem is then to show that its image lies in $R_{\nu}$. 
Similarly, the strong conjecture for $\GL_{n-1}(F)$
gives us maps $Z_{\nu'} \rightarrow R_{\nu'}$ for each factor $Z_{\nu'}$ of the integral Bernstein center of $\GL_{n-1}(F)$,
whose image we have control over.  Since these maps are compatible with local Langlands, they take the $\gamma$-factor of the pair $W_{\nu} \times W_{\nu'}$
to a Laurent series with coefficients in $R_{\nu}[\frac{1}{\ell}] \otimes R_{\nu'}$ that specializes, at any characteristic zero point $x$, to the $\gamma$-
factor of the specialization of $\rho_{\nu} \otimes \rho_{\nu'}$ at $x$.

The second key ingredient is the work of~\cite{galoisgamma}, where we show that given a family $\rho$ of Galois representations over a suitable base $R$, 
there is a unique Laurent series with coefficients in $R$ that interpolates the $\gamma$-factors of $\rho$ at characteristic zero points.  This means
that the image of the $\gamma$-factor of $W_{\nu} \times W_{\nu'}$ actually lies in $R_{\nu} \otimes R_{\nu'}$.  In particular, if $Z'$ denotes the preimage
of $R_{\nu}$ under the map $Z_{\nu'} \rightarrow R_{\nu}[\frac{1}{\ell}]$, then the coefficients of every universal $\gamma$-factor lie in $Z' \otimes R_{\nu'}$.
Corollary~\ref{cor:main descent} then shows that $Z' = Z_{\nu}$, completing the proof.

{\em Acknowledgements:}
The first author is grateful to Richard Taylor, Marie-France Vign\'{e}ras, and Jean-Fran\c{c}ois Dat for their ideas and encouragement.  The first author was partially
supported by EPSRC grant EP/M029719/1. The second author is grateful to Matthew Emerton, George Boxer, and Claus Sorensen for their interest and encouragement. 

\section{co-Whittaker modules} \label{sec:co-whittaker}

The appropriate context for the study of converse theorems in modular
representation theory was developed by the first author in
~\cite{bernstein3}.  We begin by summarizing the relevant results.

Let $F$ be a $p$-adic field, and let $G_n$ denote the group $\GL_n(F)$.
Let $k$ be an algebraically closed field of characteristic $\ell$ different
from $p$, and let $A$ be a $W(k)$-algebra.  Since $k$ contains the $p$-power
roots of unity, we can fix an additive character $\psi: F \rightarrow W(k)^{\times}$.
We will also regard $\psi$ as a character of the subgroup $U_n$ of $G_n$
consisting of unipotent upper triangular matrices, via the formula
$$\psi(u) = \psi(u_{12} + \dots + u_{n-1,n}).$$
Since $A$ is a $W(k)$-algebra, we will often regard $\psi$ as an $A$-valued character.

In this context we have an ``$n$th derivative functor'' $V \mapsto V^{(n)}$
from the category
$\Rep_A(G_n)$ of smooth $A[G_n]$-modules to the category of $A$-modules,
that takes a smooth $A[G_n]$-module $V$ to the module of $\psi$-coinvariants
in $V$.  We refer the reader to~\cite{EH-families}, section 3.1, for basic
properties of the functor $V \mapsto V^{(n)}$; in particular this functor
is exact, and we have a natural isomorphism 
$$(V \otimes_A M)^{(n)} \cong V^{(n)} \otimes_A M$$
for any $A$-module $M$.

Central to our approach is the notion of a {\em co-Whittaker $A[G_n]$-module},
defined below:

\begin{defn}[\cite{bernstein3}, 6.1] A smooth $A[G_n]$-module $V$ is co-Whittaker if
the following conditions hold:
\begin{enumerate}
\item $V$ is admissible as an $A[G_n]$-module,
\item $V^{(n)}$ is a free $A$-module of rank one, and
\item if $W$ is a quotient of $V$ such that $W^{(n)} = 0$, then $W = 0$.
\end{enumerate}
\end{defn}

(This is not quite the definition given in~\cite{bernstein3}, but is easily seen to
be equivalent, for instance by using Lemma 3.4 of~\cite{bernstein3}.)

Conditions (1)-(3) above imply easily that for any co-Whittaker $A[G_n]$-module $V$,
the map $A \rightarrow \End_{A[G_n]}(V)$ is an isomorphism (c.f.~\cite{bernstein3},
Proposition 6.2).

If $V$ and $V'$ are co-Whittaker $A[G_n]$-modules, we say that $V$ {\em dominates}
$V'$ if there is a surjection $V \rightarrow V'$ that induces an isomorphism
of $V^{(n)}$ with $(V')^{(n)}$.  This induces an equivalence relation
on the set of isomorphism classes of co-Whittaker $A[G_n]$-modules,
in which $V$ and $V'$ are equivalent if there exists a co-Whittaker module
$V''$ that dominates both $V$ and $V'$.

In~\cite{bernstein3} we construct ``universal'' co-Whittaker modules up
to this notion of equivalence.  The key tool is the integral Bernstein center
for $\GL_n(F)$; that is, the center of the category $\Rep_{W(k)}(G_n)$.

Recall that for an abelian category ${\mathcal A}$, the {\em center} of ${\mathcal A}$
is the ring of natural transformations from the identity functor on ${\mathcal A}$
to itself.  This ring acts naturally on every object in ${\mathcal A}$.  We
denote by $Z_n$ the center of $\Rep_{W(k)}(G_n)$.

A primitive idempotent $e$ of $Z_n$ gives rise to a direct factor $e\Rep_{W(k)}(G_n)$;
this is the full subcategory of $\Rep_{W(k)}(G_n)$ on which $e$ acts as the identity.
The primitive idempotents of $Z_n$ were described in~\cite{bernstein1}; they are
in bijection with inertial equivalence class of pairs $(L,\pi)$, where $L$ is a
Levi subgroup of $G_n$ and $\pi$ is an irreducible supercuspidal $k$-representation of $L$.
If $e$ is the idempotent corresponding to the pair $(L,\pi)$, then
a representation $V$ in $\Rep_{W(k)}(G_n)$ lies in $e\Rep_{W(k)}(G_n)$ if, and only if,
every simple subquotient of $V$ has {\em mod $\ell$ inertial supercuspidal support}
given by $(L,\pi)$ (in the sense of~\cite{bernstein1}, Definition 4.12).

The center $Z_n$ decomposes as a product, over the primitive idempotents $e$, of the 
rings $eZ_n$.  The structure of these rings was investigated extensively in~\cite{bernstein1};
in particular, we have:

\begin{thm}[\cite{bernstein1}, Theorem 10.8] The ring $eZ_n$ is a finitely generated,
reduced, $\ell$-torsion free $W(k)$-algebra.
\end{thm}

Let $W_n$ be the smooth $W(k)[G_n]$-module $\cInd_{U_n}^{G_n} \psi$.  Then for any
primitive idempotent $e$ of $Z_n$, we have an action of $eZ_n$ on $eW_n$.  We then have:

\begin{thm}[\cite{bernstein3}, Theorem 6.3] The smooth $eZ_n[G_n]$-module $eW_n$
is a co-Whittaker $eZ_n[G_n]$-module.
\end{thm}

Now let $A$ be a $W(k)$-algebra, and let $V$ be a co-Whittaker $A[G]$-module.  Suppose
further that $V$ lies in $e\Rep_{W(k)}(G_n)$ for some primitive idempotent $e$.  Then
$eZ_n$ acts on $V$, and since every endomorphism of $V$ is a scalar (\cite[Prop 6.2]{bernstein3}), this action
is given by a map $f_V: eZ_n \rightarrow A$.  Note that if $V$ dominates $V'$,
then the maps $f_V$ and $f_{V'}$ coincide.

In the converse direction, we have:

\begin{thm}[\cite{bernstein3}, Theorem 6.3] If $A$ is Noetherian and has an $eZ_n$-algebra structure, the module $eW_n \otimes_{eZ_n} A$
is a co-Whittaker $A[G_n]$-module that dominates $V$.
\end{thm}

In particular the maps $f_V$ and $f_{V'}$ coincide if, and only if, $V$ is equivalent
to $V'$, and in this case both $V$ and $V'$ are dominated by $eW_n \otimes_{eZ_n} A$.
We thus say that, up to the equivalence relation induced by dominance, $eW_n$ is
the universal co-Whittaker module contained in $e\Rep_{W(k)}(G_n)$.

The smooth dual of a co-Whittaker module is rarely itself a co-Whittaker module.
However, there is a natural ``duality'' operation on equivalence classes of
co-Whittaker modules.  Indeed, following Gelfand-Kazhdan, for {\em any} smooth
$W(k)[G_n]$-module $V$, we denote by $V^{\iota}$ the $W(k)[G_n]$-module with the same
underlying set as $V$, but for which the action of $g \in G_n$ on $V^{\iota}$ coincides
with the action of $(g^t)^{-1}$ on $V$.  The functor $V \mapsto V^{\iota}$ is left exact
and covariant, and, for any $V$, the $W(k)$-modules $V^{(n)}$ and $(V^{\iota})^{(n)}$ 
are isomorphic.  Moreover, $W_n^{\iota}$ is isomorphic to $W_n$.  (The latter two claims
are immediate from the fact that the character $\psi$ of $U_n$ is conjugate to the
character $u \mapsto \psi^{-1}(u^t)$ of the ``opposite'' unipotent group $\overline{U}_n$.)
In particular, if $V$ is a co-Whittaker $A[G_n]$-module, then so is $V^{\iota}$.

Moreover, Bernstein-Zelevinski show that if $V$ is an irreducible smooth $\overline{\CK}[G_n]$-module, then $V^{\iota}$ is simply the smooth dual of $V$ (this follows immediately from \cite[Thm 7.3]{BZ1}).

From this we deduce:

\begin{prop}
Let $e$ be a primitive idempotent of $Z_n$.  Then there exists a primitive
idempotent $e^{\iota}$ of $Z_n$ such that an irreducible $\overline{\CK}$ representation
of $G_n$ lies in $e\Rep_{W(k)}(G)$ if, and only if, its smooth $\overline{\CK}$-dual
lies in $e^{\iota}\Rep_{W(k)}(G)$.  Moreover, there is a unique isomorphism:
$$x \mapsto x^{\iota}: e Z_n \rightarrow e^{\iota} Z_n$$
such that for all irreducible $\overline{\CK}$-representations $\pi$ in $e\Rep_{W(k)}(G)$,
the action of $eZ_n$ on $\pi$ is given by the map: 
$$x \mapsto f_{\pi^{\vee}}(x^{\iota}): eZ_n \rightarrow \overline{\CK}.$$
(Here $\pi^{\vee}$ is the smooth dual of $\pi$.)
\end{prop}
\begin{proof}
We have seen that the map $Z_n \rightarrow \End_{W(k)[G_n]}(W_n)$ is an isomorphism.  As
$W_n$ is isomorphic to $W_n^{\iota}$, we obtain an involution $z \mapsto z^{\iota}$ on $Z_n$
by identifying $Z_n$ with $\End_{W(k)[G_n]}(W_n)$ and considering the involution
$f \mapsto f^{\iota}$ on this space of endomorphisms.  All of the claims above are now
immediate from the basic properties of the functor $V \mapsto V^{\iota}$. 
\end{proof}

\begin{prop} Let $V$ be a co-Whittaker $A[G]$-module.  Let $f_V^{\iota}: e^{\iota}Z_n \rightarrow A$
be the map defined by $f_V^{\iota}(x) = f_V(x^{\iota})$.  Then $V^{\iota}$ is isomorphic
to the co-Whittaker module $e^{\iota}W_n \otimes_{e^{\iota}Z_n, f_V^{\iota}} A$.
\end{prop}
\begin{proof}
This is clear, as $z \in Z_n$ acts on $V$ by $f_V(z)$, hence on $V^{\iota}$ by
$f_V(z)^{\iota}$.  It thus suffices to show that the latter is equal to $f_V(z^{\iota})$.
This is clear, though, by applying the functor $V \mapsto V^{\iota}$ to diagrams of the form:
$$
\begin{array}{ccc}
W_n & \rightarrow & W_n \\
\downarrow & & \downarrow\\
V & \rightarrow & V\\
\end{array}
$$
where the left-hand and right-hand maps $W_n \rightarrow V$ agree, the top map is multiplication by $z$,
and the bottom map is multiplication by $f_V(z)$.
\end{proof}

We regard the functor $V \mapsto V^{\iota}$ on co-Whittaker models as an operation that
interpolates the ``smooth dual'' operation across a co-Whittaker family.

\section{Whittaker functions and Schwartz functions} \label{sec:whittaker}

Co-Whittaker modules are useful to us because they provide a natural context
for studying the variation of Whittaker functions in families.  In this section
we recall the details of this theory.

For any smooth $A[G]$-module $V$, and any $A$-module $M$, Frobenius reciprocity
gives an isomorphism:

$$\Hom_{A}(V^{(n)},M) \cong \Hom_{A[G]}(V, \Ind_{U_n}^{G_n} \psi_M),$$
where $\psi_M$ is a copy of $M$ on which $U_n$ acts via $\psi$.  We call
an element $f$ of $\Hom_{A}(V^{(n)},M)$ an {\em $M$-valued Whittaker functional}
on $V$.  For any element $v$ of $V$, the {\em $M$-valued Whittaker function}
$w_f(v)$ attached to $f$ and $v$ is the element of $\Ind_{U_n}^{G_n} \psi$ 
given by ${\tilde f}(v)$, where ${\tilde f}$ is the element of 
$\Hom_{A[G]}(V,\Ind_{U_n}^{G_n} \psi)$ corresponding to $f$.  We will often
regard $w_f(v)$ as a smooth function on $G_n$; explicitly, one has
$w_f(v)(g) = f(\overline{gv})$, where $\overline{gv}$ is the image of $gv$
in $V^{(n)}$.

We will most often make use of this when V is a co-Whittaker module over $A$,
and $f$ is an isomorphism $V^{(n)} \cong A$.  In this case we obtain
a map ${\tilde f}: V \rightarrow \Ind_{U_n}^{G_n} \psi_A$.  We denote
the image of this map by ${\mathcal W}(V,\psi)$ and call it the
{\em $A$-valued Whittaker model} of $V$ with respect to $\psi$.

For a fixed $V$ it is possible to construct Whittaker functions with
prescribed values on a large subgroup of $G_n$.  More precisely,
let $P_n$ be the so-called ``mirabolic subgroup'' of $G_n$; that is,
the subgroup of $G_n$ consisting of matrices whose last row
is the vector $(0,\dots,0,1)$.

There are natural isomorphisms of functors:
$$V^{(n)} \cong \Hom_{W(k)[P_n]}(\cInd_{U_n}^{P_n} \psi, V)$$
$$\Hom_{W(k)}(V^{(n)},W) \cong \Hom_{W(k)[P_n]}(V,\Ind_{U_n}^{P_n} \psi_W)$$
due to Bernstein-Zelevinsky.  More precisely, this follows from the existence
of isomorphisms:
$$(\Phi^+)^{n-1}\Psi^+ W(k) \cong \cInd_{U_n}^{P_n} \psi$$
$$(\hat \Phi^+)^{n-1}\Psi^+ W(k) \cong \Ind_{U_n}^{P_n} \psi$$
(notation as in~\cite{BZ}, section 3) and the formalism of the Bernstein-Zelevinsky functors
(developed in~\cite{BZ}, section 3 over the complex numbers and in~\cite{EH-families},
section 3.1 over $W(k)$.)

The first of these isomorphisms give rise to a natural map:
$$\cInd_{U_n}^{P_n} \psi \otimes V^{(n)} \rightarrow V.$$
Note that we may identify $\cInd_{U_n}^{P_n} \psi \otimes V^{(n)}$
with $\cInd_{U_n}^{P_n} \psi_{V^{(n)}}$, allowing us to view the domain
of this map as a space of compactly supported $V^{(n)}$-valued functions on $P_n$.
The image of this map is often called the space of {\em Schwartz functions}
of $V$.

On the other hand, the identity map $V^{(n)} \rightarrow V^{(n)}$
gives rise to a map $V \rightarrow \Ind_{U_n}^{G_n} \psi_{V^{(n)}}$.  Composing with
restriction to $P_n$ gives us a series of maps:
\begin{equation}\label{eq:1}
\cInd_{U_n}^{P_n} \psi_{V^{(n)}} \rightarrow V \rightarrow \Ind_{U_n}^{G_n} \psi_{V^{(n)}}
\rightarrow \Ind_{U_n}^{P_n} \psi_{V^{(n)}}
\end{equation}
and we have:

\begin{lemma} \label{lemma:include}
The composition of the chain of maps~\ref{eq:1} is the natural inclusion:
$$\cInd_{U_n}^{P_n} \psi_{V^{(n)}} \rightarrow \Ind_{U_n}^{P_n} \psi_{V^{(n)}}.$$
\end{lemma}
\begin{proof}
This is a fairly easy consequence of the Bernstein-Zelevinski formalism.
The composition of the maps:
$$
V \rightarrow \Ind_{U_n}^{G_n} \psi_{V^{(n)}}
\rightarrow \Ind_{U_n}^{P_n} \psi_{V^{(n)}}
$$
is simply the map $V \rightarrow \Ind_{U_n}^{P_n} \psi_{V^{(n)}}$
attached, via Frobenius reciprocity, to the identity map on $V^{(n)}$.  Translating this
into the language of~\cite{BZ} via the isomorphisms $V^{(n)} \cong \Psi^-(\Phi^-)^{(n-1)} V$
and $\Ind_{U_n}^{P_n} \psi_W \cong (\hat \Phi^+)^{(n-1)} \Psi^+ W$ we find that, under these isomorphisms,
this map corresponds to the map:
$$V \rightarrow (\hat \Phi^+)^{(n-1)} \Psi^+ \Psi^- (\Phi^-)^{(n-1)} V$$
which comes from the identity on $\Psi^- (\Phi^-)^{(n-1)} V,$ via the adjunctions of the pairs 
$\Phi^-,{\hat \Phi}^+$, and $\Psi^-,\Psi^+$.

Similarly, under the isomorphism of $\cInd_{U_n}^{P_n} \psi_W$ with $(\Phi^+)^{(n-1)}\Psi^+ W$,
the map $\cInd_{U_n}^{P_n} \psi_{V^{(n)}} \rightarrow V$ is the map:
$$(\Phi^+)^{(n-1)} \Psi^+ \Psi^- (\Phi^-)^{(n-1)} V \rightarrow V$$
corresponding to the identity on $\Psi^- (\Phi^-)^{(n-1)} V$ via the adjunction of $\Phi^-$ and $\Phi^+$.

It follows that the composition of these two maps is given by the inclusion of $(\Phi^+)^{(n-1)}$ into
$({\hat \Phi}^+)^{(n-1)}$, and, under the identifications we have made, inclusion corresponds in turn to 
the inclusion of $\cInd_{U_n}^{P_n} \psi_{V^{(n)}}$ in $\Ind_{U_n}^{P_n} \psi_{V^{(n)}}$.
\end{proof}

It is immediate from Lemma~\ref{lemma:include} that when $V$ is a co-Whittaker $A[G]$-module, 
and we fix an isomorphism $f$ of $V^{(n)}$ with $A$,
then given any function $h$ in $\cInd_{U_n}^{P_n} \psi_A$, there exists an element $v$
of $V$ such that the restriction of $w_f(v)$ to $P_n$ is equal to $h$.

Our main interest in the Schwartz functions and Whittaker functions comes from the
following ``descent'' result for co-Whittaker modules:

\begin{thm} \label{thm:descent}
Let $V$ be a co-Whittaker $A[G]$-module in $e\Rep_{W(k)}(G)$, and fix an 
isomorphism $f: V^{(n)} \cong A$.
Let $A'$ be a $W(k)$-subalgebra of $A$, and suppose that the composition:
$$\cInd_{U_n}^{P_n} \psi_{A'} \rightarrow \cInd_{U_n}^{P_n} \psi_A \rightarrow V
\rightarrow \Ind_{U_n}^{G_n} \psi_A$$
has image contained in $\Ind_{U_n}^{P_n} \psi_{A'}$ (i.e. the ``$A'$-valued
Schwartz functions of $V$ have $A'$-valued Whittaker functions'').
Then the map $eZ_n \rightarrow A$ giving the action of $eZ_n$ on $V$ factors through
$A'$.
\end{thm}
\begin{proof}
Let $V'$ be the preimage of $\Ind_{U_n}^{G_n} \psi_{A'}$ under the map
$$V \rightarrow \Ind_{U_n}^{G_n} \psi_A.$$
The map $V' \rightarrow \Ind_{U_n}^{G_n} \psi_{A'}$ gives rise to
a map $(V')^{(n)} \rightarrow A'$ that fits into a commutative diagram:
$$
\begin{array}{ccc}
V'^{(n)} & \rightarrow & A'\\
\downarrow & & \downarrow\\
V^{(n)} & \rightarrow & A
\end{array}
$$
in which the left-hand vertical map is injective (by exactness of the derivative)
and the bottom horizontal map is the isomorphism $f$.  We therefore
have an injection of $(V')^{(n)}$ into $A'$.

On the other hand, our assumptions above show that the image
of $\cInd_{U_n}^{P_n} \psi_{A'}$ in $V$ is contained in $V'$; the
map 
$$\cInd_{U_n}^{P_n} \psi_{A'} \rightarrow V'$$
gives rise to a map $A' \rightarrow (V')^{(n)}$ whose composition
with the injection of $(V')^{(n)}$ into $A'$ is the identity.
The isomorphism $V^{(n)} \rightarrow A$ thus identifies $(V')^{(n)}$
with $A'$.

The ring $eZ_n$ acts $A'$-linearly on $V'$ and hence on $(V')^{(n)}$; this
action is given by a map $eZ_n \rightarrow A'$.  Since $(V')^{(n)}$ generates
$V^{(n)}$ as an $A$-module, and the map $(V')^{(n)} \rightarrow V^{(n)}$
is $eZ_n$-equivariant, the action of $eZ_n$ on $V^{(n)}$ is via this same map, and
the result is now immediate.
\end{proof}

In light of this result it will be useful to develop criteria for when
a Whittaker function of $V$ takes values in a given subalgebra $A'$ of $A$.
We first note:

\begin{lemma} \label{lem:integrate}
Let $A'$ be a $W(k)$-subalgebra of $A$, and let $h$ be a function in the compact induction $\cInd_{U_n}^{G_n} \psi_A$.
Suppose that for all primitive idempotents $e$ of $Z_n$, and
all functions $g$ in $e\cInd_{U_n}^{G_n} \psi^{-1}$, the integral
$$\int_{U_n\setminus G_n} h(x)g(x) d\mu(x)$$
takes values in $A'$.  Then $h$ lies in $\cInd_{U_n}^{G_n} \psi_{A'}$.
\end{lemma}
\begin{proof}
As $\cInd_{U_n}^{G_n} \psi^{-1}$ is the direct sum, over all primitive $e$, of $e\cInd_{U_n}^{G_n} \psi^{-1}$,
it suffices to show, for any $h$ that does {\em not} take values in $A'$, there exists a
$g \in \cInd_{U_n}^{G_n} \psi^{-1}$ such that the integral of $h(x)g(x)$ does not lie in $A'$.
Fix such an $h$, and
a compact open subgroup $K$ of $G_n$ such that $h$ is fixed under right translation by all $k \in K$.
Then (since the values of $\psi$ lie in $A'$) there is an element $x$ of $G_n$ such that
$h$ takes values in $A \setminus A'$ on the double coset $U_n x K$.  On the other hand, for $K$ sufficiently
small there exists $g: U_n x K \rightarrow A'$ such that $g(uxk) = \psi^{-1}(u)$ for all $u \in U_n$, $k \in K$.
It is then clear that the integral of $h(x)g(x)$ is valued in $A \setminus A'$.
\end{proof}

We will use the fact that the spaces $e\cInd_{U_n}^{G_n} \psi^{-1}$ are co-Whittaker
modules to convert this into a criterion that involves integration against Whittaker functions.
Let $A$ and $B$ be $W(k)$-algebras, let $h$ be an element of $\cInd_{U_n}^{G_n} \psi_A$, and
let $V$ be a co-Whittaker $B[G_n]$-module.  Fixing an isomorphism $f: V^{(n)} \cong B$
and an element $v$ of $V$ gives us a Whittaker function $w_f(v)$ in $\Ind_{U_n}^{G_n} \psi_B$. Given any Whittaker function $w$ in ${\mathcal W}(V,\psi)$ we can define the dual Whittaker function ${\tilde w}$ by ${\tilde w}(g) = w(\delta_n(g^t)^{-1})$, where $\delta_n$ is the matrix
in $G_n$ whose only nonzero entries are $1$'s along the antidiagonal. We have:
\begin{prop} 
Let $w$ lie in ${\mathcal W}(V,\psi)$.  Then
$\tilde w$ lies in ${\mathcal W}(V^{\iota},\psi^{-1})$.
\end{prop}
\begin{proof}
This is clear, as $G_n$ acts on Whittaker functions by right translation, so one has $\tilde{gw} = (g^t)^{-1}\tilde{w}$ for all $g$.
\end{proof}

Let $\tilde{w}_f(v)$ be the dual Whittaker function to $w_f(v)$. Then the integral
$$\int_{U_n\setminus G_n} h(x) \otimes \tilde{w}_f(v)(x) d\mu(x)$$
makes sense as an element of $A \otimes B$.

This is particularly useful in the case where $B = eZ_n$ and $V$ is the universal
co-Whittaker module $eW_n$.  In this case one has some useful additional structure.
First, regarding elements of $W_n = \cInd_{U_n}^{G_n} \psi$, we see that
``evaluation at the identity'' gives a canonical map $V^{(n)} \rightarrow W(k)$
that corresponds to the inclusion of $\cInd_{U_n}^{G_n} \psi$ in $\Ind_{U_n}^{G_n} \psi$.
We denote this map by $\theta_{e,n}$.

We also note that there is a natural isomorphism: $\End_{W(k)[G]}(eW_n) \cong eW_n^{(n)}$.
Since $eW_n$ is a co-Whittaker $eZ_n[G]$-module, its endomorphisms
are canonically isomorphic to $eZ_n$, and we thus canonically identify
$$eZ_n \cong \End_{W(k)[G]}(eW_n) \cong eW_n^{(n)}.$$
In particular we can (and do) regard $\theta_{e,n}$ as a $W(k)$-linear map $eZ_n \rightarrow W(k)$.

Let $f$ be the above isomorphism $(eW_n)^{(n)}\cong eZ_n$. Immediately from the definitions one finds that for $v \in eW_n$ and $x \in G$, one has $v(x) = \theta_{e,n}(w_f(v)(x))$;
that is, one can recover the function $v: G_n \rightarrow W(k)$ from the
$eZ_n$-valued Whittaker function of $v$.  Similarly, one has $\tilde{v}(x) = \theta_{e,n}(\tilde{w}_f(v)(x))$, where $\tilde{v}(x)=v(\delta_n(g^t)^{-1})$.

From this, together with Lemma~\ref{lem:integrate}, we immediately deduce:
\begin{prop} \label{prop:integrate}
Let $A'$ be a $W(k)$-subalgebra of $A$, and let $h$ be a function in $\cInd_{U_n}^{G_n} \psi_A$.
Suppose that for all primitive idempotents $e$ of $Z_n$, and all functions $g$ in $e\cInd_{U_n}^{G_n} \psi$,
the integral
$$I(h,g) = \int_{U_n\setminus G_n} h(x) \otimes \tilde{w}_f(g) (x) d\mu(x)$$
(an element of $A \otimes eZ_n$)
satisfies $(1_A \otimes \theta_{e,n})(I(h,g)) \in A'$.  Then $h$ lies in $\cInd_{U_n}^{G_n} \psi_{A'}$.
\end{prop}

Another useful observation, also immediate from the above, is:
\begin{cor} \label{cor:integrate}
Let $A'$ be a $W(k)$-subalgebra of $A$, and let $h$ be a function in $\cInd_{U_n}^{G_n} \psi_A$.
Suppoase that for all primitive idempotents $e$ of $Z_n$, and all functions $g$ in $e\cInd_{U_n}^{G_n} \psi$,
the integral $I(h,g)$ defined above lies in $A' \otimes eZ_n$.  Then $h$ lies in $\cInd_{U_n}^{G_n} \psi_{A'}$.
\end{cor}

\section{Zeta integrals and $\gamma$-factors} \label{sec:gamma}

We now recall recent work of the second author~\cite{pairs} that constructs
zeta integrals and $\gamma$-factors attached to pairs of co-Whittaker modules.
Let $A_1$ and $A_2$ be Noetherian $W(k)$-algebras, and let $V_1$, $V_2$ be a co-Whittaker $A_1[G_n]$-module,
and a co-Whittaker $A_2[G_{n-1}]$-module, respectively.  Let $R$ denote the $W(k)$-algebra 
$A_1 \otimes A_2$.  Fixing isomorphisms of $V_1^{(n)}$ with $A_1$ and $V_2^{(n)}$ with $A_2$
gives us Whittaker functionals on $V_1$ and $V_2$.

For elements $w_1$ of ${\mathcal W}(V_1,\psi)$ and $w_2$ of ${\mathcal W}(V_2,\psi^{-1})$ we define
the ``zeta integral'' $\Psi(w_1,w_2,X)$ to be the formal series $\sum_{m = -\infty}^{\infty} c_m X^m,$
where the coefficient $c_m$ is given by the integral:
$$c_m := \int_{U_{n-1} \setminus \{g \in G_{n-1}: v(\det g) = m\}} 
\left(w_1(\begin{smallmatrix} g & 0\\0 & 1\end{smallmatrix}) \otimes w_2(g)\right) dg.$$

By~\cite{pairs}, Lemma 2.1, $\Psi(w_1,w_2,X)$ is a well-defined element of $R[[X]][X^{-1}]$.  Note that the formation
of zeta integrals is compatible with base change, in the following sense: if we have a map $f:A_1 \rightarrow B_1$,
then the function $f(w_1)$ defined by $f(w_1)(x) = f(w_1(x))$ is an element of ${\mathcal W}(V \otimes_{A_1,f} B_1,\psi)$,
and if $\Psi(w_1,w_2,X)$ is given by $\sum c_i X^i$, then $\Psi(f(w_1),w_2,X)$ is given by $\sum (f \otimes 1)(c_i) X^i$.
A similar statement holds for base change via maps $A_2 \rightarrow B_2$.

A key point of~\cite{pairs} is a rationality result for such zeta integrals.  Let $S$ be the multiplicative system
in $R[X,X^{-1}]$ consisting of all polynomials whose leading and trailing coefficients are units.  Such polynomials
are units in $R[[X]][X^{-1}]$, and nonzerodivisors in $R[X,X^{-1}]$.  We therefore obtain an embedding:
$$S^{-1}R[X,X^{-1}] \rightarrow R[[X]][X^{-1}].$$

When $R$ is a Noetherian ring, the following lemma characterizes the image of $S^{-1}R[X,X^{-1}]$ in $R[[X]][X^{-1}]$:

\begin{lemma}
Suppose $R$ is Noetherian.
Let $f$ be an element of $R[[X]][X^{-1}]$, and let $M_f$ be the $R[X,X^{-1}]$-submodule of
$R[[X]][X^{-1}]/R[X,X^{-1}]$ generated by $f$.  Then $f$ lies in $S^{-1}R[X,X^{-1}]$ if, and only if,
$M_f$ is a finitely generated $R$-module.
\end{lemma}
\begin{proof}
Suppose that $M_f$ is a finitely generated $R$-module.  Then there exists a positive integer $d$ such that
both $X^{-1}f$ and $X^m f$ are in the $R$-submodule of $M_f$ generated by $f,Xf, \dots, X^{d-1}f$.  We may
thus write $X^{-1}f = P(X) f$ and $X^d f = Q(X)f$ where $P$ and $Q$ lie in $R[X]$ and have degree at most $d-1$.
Then $X^{-1} - P(X) - Q(X) + X^d$ annihilates $f$ in $M_f$ and lies in $S$.  Thus $f$ lies in $S^{-1}R[X,X^{-1}]$

Conversely, suppose $f$ lies in $S^{-1}R[X,X^{-1}]$.  Then there is an element of $S$ that annihilates $f$ in $M_f$,
which we may take, without loss of generality, to be a polynomial of degree $d$.  Then $M_f$ is spanned over $R$
by $f, Xf, \dots, X^{d-1}f$.
\end{proof}

\begin{cor} \label{cor:finite}
Let $R'$ be a Noetherian $W(k)$-subalgebra of $R$ such that $R$ is finitely generated as an $R'$-module.
Let $S'$ be the subset of $R'[X,X^{-1}]$ consisting of polynomials whose first and last nonzero coefficients are units in $R'$.
Then $(S')^{-1}R'[X,X^{-1}]$ is the intersection, in $R[[X]][X^{-1}]$, of the subrings $R'[[X]][X^{-1}]$ and $S^{-1}R[X,X^{-1}]$.
\end{cor}
\begin{proof}
It is clear that $(S')^{-1}R'[X,X^{-1}]$ is contained in this intersection.  Conversely, suppose $f$ is in this
intersection.  Then $M_f$ is finitely generated over $R$, and hence finitely generated over $R'$.  But then
$f$ lies in $(S')^{-1}R'[X,X^{-1}]$ as claimed.
\end{proof}

The rationality result of~\cite{pairs} now states:
\begin{thm}[\cite{pairs}, Theorem 2.2]
For any $w_1 \in {\mathcal W}(V_1,\psi)$ and $w_2 \in {\mathcal W}(V_2,\psi^{-1})$, the zeta integral
$\Psi(w_1,w_2,X)$ lies in $S^{-1}R[X,X^{-1}]$.
\end{thm}

The zeta integrals attached to pairs of co-Whittaker modules satisfy a functional equation generalizing the
functional equation over the complex numbers:
\begin{thm}[\cite{pairs}, Theorem 3.4]
There exists a unique element $\gamma(V_1 \times V_2, X, \psi)$ of $S^{-1}R[X,X^{-1}]$ such that, for all $w_1$ in ${\mathcal W}(V_1,\psi)$
and $w_2$ in ${\mathcal W}(V_2,\psi^{-1})$, one has:
$$\Psi(w_1,w_2,X)\gamma(V_1 \times V_2, X, \psi) \omega_{V_2}(-1)^{(n-1)} = \Psi({\tilde w_1}, {\tilde w_2}, X^{-1}).$$
(Here $\omega_{V_2}$ is the central character of $V_2$, viewed as taking values in $A_2$.) 
\end{thm}

Note that the compatibility of zeta integrals with base change, together with uniqueness of $\gamma$-factors, immediately implies
a similar compatbility of $\gamma$-factors with base change.  

\section{A descent result} \label{sec:descent}

Now fix a Noetherian $W(k)$-algebra $A$, and a co-Whittaker $A[G_n]$-module $V$ in $e\Rep_{W(k)}(G_n)$.  For each
idempotent $e'$ of $Z_{n-1}$, the module $e'W_{n-1}$ is a co-Whitaker $e'Z_{n-1}[G_{n-1}]$-module,
and one can form the rational function $\gamma(V \times e'W_{n-1}, X^{-1}, \psi)$.  If we expand this as a
power series in $X$, then the coefficients of this power series
lie in $A \otimes e'Z_{n-1}$.  

Compatibility of $\gamma$-factors with base change implies that if $V$ arises by base change from some subalgebra $A'$ of
$A$, then the coefficients of $\gamma(V \times e'W_{n-1}, X^{-1}, \psi)$ lie in $A' \otimes e'Z_{n-1}$,
as do the coefficients of $\gamma(V^{\iota} \times e'W_{n-1}, X, \psi)$.
The objective of this section
is to prove a partial converse to this theorem.

Fix a $W(k)$-subalgebra $A'$ of $A$ such that $A$ is finitely generated as an $A'$-module, and let $V$
be a co-Whittaker $A[G_n]$-module such that for all $e'$, the coefficients of $\gamma(V \times e'W_{n-1}, X^{-1}, \psi)$
and $\gamma(V^{\iota} \times e'W_{n-1}, X, \psi)$
lie in $A' \otimes eZ_{n-1}$.  Let $R$ denote the ring $(A \otimes e'Z_{n-1})$, and let $R'$ denote
the subalgebra $(A' \otimes e'Z_{n-1})$.  As usual $S$ and $S'$ will denote the multiplicative subsets
of $R[X,X^{-1}]$ and $R'[X,X^{-1}]$ consisting of polynomials whose first and last coefficients are units.

We then have:
\begin{prop} \label{prop:tilde}
Let $w$ be an element of ${\mathcal W}(V,\psi)$.  Then
the restriction of $w$ to $P_n$ takes values in $A'$ if and only if the restriction of ${\tilde w}$
to $P_n$ takes values in $A'$.
\end{prop}
\begin{proof}
Fix a primitive idempotent $e'$, and an element $w_2$ of ${\mathcal W}(e'W_{n-1},\psi^{-1})$.
The functional equation for $\Psi(w,{\tilde w_2}, X)$ reads:
$$\Psi(w,{\tilde w_2},X) \gamma(V \times (e')^{\iota}W_{n-1},X,\psi) = \pm \Psi({\tilde w_1}, w_2, X^{-1}).$$
Applying the involution $X \mapsto X^{-1}$ of $S^{-1}R[X,X^{-1}]$ we obtain the identity:
$$\Psi(w,{\tilde w_2},X^{-1}) \gamma(V \times (e')^{\iota}W_{n-1},X^{-1},\psi) = \pm \Psi({\tilde w}, w_2, X).$$

Since the zeta integral $\Psi(w,{\tilde w_2},X)$ depends, by definition, only on the restriction of $w$
to $P_n$, this integral lies in $R'[[X]][X^{-1}]$.  Since it also lies in $S^{-1}R[X,X^{-1}]$, and $R$ is finitely
generated over $R'$, we find that $\Psi(w,{\tilde w_2},X)$ lies in the subring $(S')^{-1}R'[X,X^{-1}]$, and hence
so does $\Psi(w,{\tilde w_2},X^{-1})$.  The latter, considered as a power series in $X$, thus lies in
$R'[[X]][X^{-1}]$.  Therefore, so does $\Psi({\tilde w}, w_2, X)$.

On the other hand, returning to the definition of the zeta integral, one sees easily that one can write
$$\int_{U_{n-1}\setminus G_{n-1}} {\tilde w}|_{G_{n-1}}(x) \otimes w_2(x) d\mu(x)$$
as a sum of coefficients of $\Psi({\tilde w}, w_2, X,\psi)$.  It follows that the integral
$$\int_{U_{n-1}\setminus G_{n-1}} {\tilde w}|_{G_{n-1}}(x) \otimes w_2(x) d\mu(x)$$
lies in $R'$.  Since this is true for all $e'$ and all $w_2$, we conclude that ${\tilde w}$ takes values in $A'$
when restricted to $G_{n-1}$, and hence also to $P_n$.

The converse argument is nearly identical, starting with the functional equation
$$\Psi({\tilde w},{\tilde w_2}, X) \gamma(V^{\iota} \times (e')^{\iota} W_{n-1},X,\psi) = \pm \Psi(w_1,w_2,X^{-1}).$$
\end{proof}

Our main results are immediate consequences of this proposition:
\begin{cor}
If $w$ is an element of ${\mathcal W}(V,\psi)$ whose restriction
to $P_n$ takes values in $A'$, then $w$ takes values in $A'$ on all of $G_n$.
\end{cor}
\begin{proof}
The set of $w$ whose restriction to $P_n$ takes values in $A'$ is clearly stable under $P_n$, and similarly
the set of $w$ such that the restriction of $\tilde w$ to $P_n$ takes values in $A'$ is stable under the transposed
group $P_n^t$.  Since $P_n$ and $P_n^t$ together generate all of $G_n$, the set of $w$ whose restriction to $P_n$
takes values in $A'$ is stable under $G_n$.  Hence any such $w$ takes values in $A'$ on all of $G_n$.
\end{proof}

\begin{cor} \label{cor:main descent}
The map $f_V: eZ_n \rightarrow A$ giving the action of $eZ_n$ on $V$ factors through $A'$.
\end{cor}
\begin{proof}
This is immediate from the previous corollary and Theorem~\ref{thm:descent}.
\end{proof}

\section{Generators for $eZ_n$} \label{sec:generators}
Although it will not be necessary for our main results, it is an interesting question, given a co-Whittaker
$A[G_n]$-module $V$, to determine the image of the map $eZ_n \rightarrow A$ giving the action of the Bernstein
center on $V$ in terms of the $\gamma$-factors of $V$.  One can try to approach this question via the techniques of
the previous section, but the finiteness hypotheses that appear in places (particularly Corollary~\ref{cor:finite})
prevent us from being completely successful.  We outline the basic ideas here:

The key approach is to construct elements of $A$ that must lie in the image of the map $eZ_n \rightarrow A$.  For instance, if $V$ arises by base change from $A'$, then $A'$ contains $(1 \otimes \theta_{e',n-1})(zc_i)$ for all $e'$, all $z\in e'Z_{n-1}$, and all coefficients $c_i$ of the power series $\gamma(V \times e'W_{n-1},X,\psi)$. Similarly,
if $d_i$ are the cofficients of $\gamma(V^{\iota} \times e'W_{n-1}, X^{-1}, \psi)$, then
$A'$ contains $(1 \otimes \theta_{e',n-1})(zd_i)$ for all $e'$, $z$, and $d_i$.  

Therefore, let $A'$ be a $W(k)$-subalgebra of $A$ that contains $(1 \otimes \theta_{e',n-1})(zc_i)$ and $(1 \otimes \theta_{e',n-1})(zd_i)$
for all choices of $e'$, $z$, and $i$.  Note that for all $e'$, $e'Z_{n-1}$ is flat over $W(k)$, so that $A' \otimes e'Z_{n-1}$ is a subalgebra
of $A \otimes e'Z_{n-1}$.

\begin{lemma}
Fix a primitive idempotent $e'$ of $Z_{n-1}$, and let $R = A \otimes e'Z_{n-1}$, $R' = A' \otimes e'Z_{n-1}$.  Let
$\Psi$ be an element of $R'[[X]][X^{-1}]$.  Then for any coefficient $c$ of $\Psi \gamma(V \times e'W_{n-1}, X^{-1}, \psi)$ or
of $\Psi \gamma(V^{\iota} \times e'W_{n-1}, X^{-1}, \psi)$, the element 
$(1 \otimes \theta_{e',n-1})(c)$ of $A$ lies in $A'$.
\end{lemma}
\begin{proof}
The coefficient $c$ is a sum of products of a cofficient of $\phi$ and a coefficient of $\gamma(V \times e'W_{n-1},X,\psi)$,
and therefore an $R'$-linear combination of coefficients of $\gamma$.  Thus $c$ is an $A'$-linear combination of elements of
$R$ of the form $zc_i$, where $z \in e'Z_{n-1}$ and $c_i$ is a coefficient of $\gamma(V \times e'W_{n-1},X,\psi)$.  It follows
(by $A$-linearity of $1 \otimes \theta_{e',n-1}$) that $(1 \otimes \theta_{e',n-1})(c)$ is an element of $A'$.
\end{proof}

Then as a slight variation on Proposition~\ref{prop:tilde} above, we obtain:
\begin{prop} \label{prop:second tilde}
Suppose that $A$ is a finitely generated $A'$-module, and let $w$ be an element of ${\mathcal W}(V,\psi)$.  Then
the restriction of $w$ to $P_n$ takes values in $A'$ if, and only if, the restriction of ${\tilde w}$ to $P_n$
takes values in $A'$.
\end{prop}
\begin{proof}
Fix a primitive idempotent $e'$, and an element $w_2$ of $W(e'W_{n-1},\psi^{-1})$.
The functional equation for $\Psi(w,{\tilde w_2},X)$ reads:
$$\Psi(w,{\tilde w_2},X) \gamma(V \times (e')^{\iota}W_{n-1},X,\psi) = \pm \Psi({\tilde w_1}, w_2, X^{-1}).$$
Applying the involution $X \mapsto X^{-1}$ of $S^{-1}R[X,X^{-1}]$ we obtain the identity:
$$\Psi(w,{\tilde w_2},X^{-1}) \gamma(V \times (e')^{\iota}W_{n-1},X^{-1},\psi) = \pm \Psi({\tilde w}, w_2, X).$$

Since the zeta integral $\Psi(w,{\tilde w_2},X)$ depends, by definition, only on the restriction of $w$
to $P_n$, this integral lies in $R'[[X]][X^{-1}]$.  Since it also lies in $S^{-1}R[X,X^{-1}]$, and $R$ is finitely
generated over $R'$, we find that $\Psi(w,{\tilde w_2},X)$ lies in the subring $(S')^{-1}R'[X,X^{-1}]$, and hence
so does $\Psi(w,{\tilde w_2},X^{-1})$.  The latter, considered as a power series in $X$, thus lies in 
$R'[[X]][X^{-1}]$.  Hence, by the lemma above, if $c$ is any coefficient of $\Psi({\tilde w}, w_2, X)$,
then $(1 \otimes \theta_{e',n-1})(c)$ lies in $A'$.

On the other hand, returning to the definition of the zeta integral, one sees easily that one can write
$$\int_{U_{n-1}\setminus G_{n-1}} {\tilde w}|_{G_{n-1}}(x) \otimes w_2(x) d\mu(x)$$
as a sum of coefficients of $\Psi({\tilde w}, w_2, X)$.  It follows that the integral
$$\int_{U_{n-1}\setminus G_{n-1}} (1 \otimes \theta_{e',n-1}){\tilde w}|_{G_{n-1}}(x) \otimes w_2(x) d\mu(x)$$
lies in $A'$.  Since this is true for all $e'$ and all $w_2$, we conclude that ${\tilde w}$ takes values in $A'$
when restricted to $G_{n-1}$, and hence also to $P_n$.

As before, the converse argument is nearly identical, starting with the functional equation
$$\Psi({\tilde w},{\tilde w_2}, X) \gamma(V^{\iota} \times (e')^{\iota} W_{n-1},X,\psi) = \pm \Psi(w_1,w_2,X^{-1}).$$
\end{proof}

Just as in the previous section, we now find:
\begin{cor}
If $A$ is a finitely generated $A'$-module, and $w$ is an element of ${\mathcal W}(V,\psi)$ whose restriction
to $P_n$ takes values in $A'$, then $w$ takes values in $A'$ on all of $G_n$.
\end{cor}

\begin{cor}
Suppose that $A$ is a finitely generated $A'$-module, or that $A$ is a complete Noetherian local ring with residue field
$k$, and $A'$ is closed in $A$.  Then the map $f_V: eZ_n \rightarrow A$ giving the action of $eZ_n$ on $V$ factors through $A'$.
\end{cor}
\begin{proof}
When $A$ is a finitely generated $A'$-module, this is immediate from the previous corollary and Theorem~\ref{thm:descent}.
When $A$ is complete local with maximal ideal ${\mathfrak m}$, set $A_r = A/{\mathfrak m}^r$, let $A'_r$ be the image
of $A'$ in $A_r$, and set $V_r = V \otimes_A A_r$.  Then $A_r$ has finite length over $W(k)$, and in particular
is a finitely generated $A'_r$-module.  Thus, for all $r$, the map $f_{V_r}$ factors through $A'_r$.  Passing to the
limit (and using the fact that $A'$ is closed in $A$), we find that $f_V$ factors through $A'$.
\end{proof}

Now let ${\mathfrak m}$ be a maximal ideal of $eZ_n$ with residue field $k$.  We then have:
\begin{cor} \label{cor:generators}
The completion of $eZ_n$ at ${\mathfrak m}$ is toplogically generated by elements of the form
$(1 \otimes \theta_{e',n-1})(zc_i)$ and $(1 \otimes \theta_{e',n-1})(zd_i)$, where $e'$ runs over the primitive
idempotents of $Z_{n-1}$, $z$ runs over the elements of $e'Z_{n-1}$, and $c_i$ and $d_i$ are the coefficients
of the power series $\gamma(eW_n \times e'W_{n-1},X^{-1})$ and $\gamma(e^{\iota}W_n \times e'W_{n-1},X^{-1})$, respectively.
\end{cor}
\begin{proof}
Apply the previous corollary with $A$ equal to the completion of $eZ_n$ at ${\mathfrak m}$, 
$V = eW_n \otimes_{eZ_n} A$, and $A'$ the closed subalgebra of $A$ generated by the elements
listed above.  Then the map $eZ_n \rightarrow A$ factors through $A'$; since $eZ_n$ is dense in $A$ the result follows.
\end{proof}

\section{The local Langlands correspondence in families} \label{sec:langlands}

Now that we have constructed a set of generators for $eZ_n$, we turn to other applications.  In particular, we
recall Conjecture 1.3.1 of~\cite{EH-families} (or rather, a slight reformulation of it in the spirit of section 7
of~\cite{bernstein3}):

\begin{conj} Let $A$ be a reduced complete Noetherian local $W(k)$-algebra, flat over $W(k)$, with residue field
$k$, and let $\rho: G_F \rightarrow \GL_n(A)$ be a continuous $n$-dimensional representation of the Galois group $G_F$.
Then there exists a (necessarily unique) admissible $W(k)[G_n]$-module $\pi(\rho)$ such that:
\begin{enumerate}
\item $\pi(\rho)$ is $A$-torsion free,
\item $\pi(\rho)$ is a co-Whittaker $A[G]$-module, and
\item for each minimal prime ${\mathfrak a}$ of $A$, the representation $\pi(\rho)_{\mathfrak a}$ is $\kappa({\mathfrak a})$-dual
to the representation that corresponds to $\rho_{\mathfrak a}^{\vee}$ under the Breuil-Schneider generic local Langlands
correspondence.
\end{enumerate}
(Here $\kappa({\mathfrak a})$ denotes the field of fractions of $A/{\mathfrak a}$, and $\rho_{\mathfrak a}^{\vee}$ is the
$\kappa({\mathfrak a})$-dual of $\rho_{\mathfrak a}$.)
\end{conj}

The main result of~\cite{bernstein3} (Theorem 7.9) gives a reformulation of this conjecture in terms of the Bernstein
center.  More precisely, fix $\rho$ as in the conjecture and let $\rhobar: G_F \rightarrow \GL_n(k)$ be its residual
representation.  Let $R_{\rhobar}^{\Box}$ denote the universal framed deformation ring of $\rhobar$; in particular
there is a map $\rho^{\Box}: G_F \rightarrow \GL_n(R_{\rhobar}^{\Box})$ lifting $\rhobar$ that is universal
for lifts of $\rhobar$ to complete Noetherian local $W(k)$-algebras with residue field $k$.

Then Conjecture 7.5 of~\cite{bernstein3} asserts that there exists a map $\LL_{\rhobar}: Z_n \rightarrow R_{\rhobar}^{\Box}$
that is ``compatible with local Langlands'', i.e., that for any map $x: R_{\rhobar}^{\Box} \rightarrow \overline{\CK}$,
the composition $x \circ \LL_{\rhobar}$ is the map $Z_n \rightarrow \overline{\CK}$ giving the action of $Z_n$ on $\Pi_x$,
where $\Pi_x$ corresponds to $\rho_x$ under local Langlands.  By Theorem 7.9 of~\cite{bernstein3}, this conjecture
implies Conjecture 1.3.1 of~\cite{EH-families}.

The paper~\cite{curtis} gives a further refinement of these conjectures.  In particular, if one fixes $\rhobar$
(or even merely the restriction $\nubar$ of $\rhobar$ to prime-to-$\ell$ inertia) then this determines a primitive idempotent
$e$ of $Z_n$.  This primitive idempotent is characterized by the following property: let $\Pi$ be an irreducible
representation of $G_n$ over $\overline{\CK}$, and let $\rho$ correspond to $\Pi$ via local Langlands.  Let $\nu$
denote the unique lift of $\nubar$ to a representation over $\overline{\CK}$.
Then $\Pi$ lies in the block corresponding to $e$ if, and only if,
the restriction of $\rho$ to the prime-to-$\ell$ inertia group $I_F^{(\ell)}$ of $F$ is isomorphic to $\nu$.

Section 9 of~\cite{curtis} constructs a finitely generated, reduced, $\ell$-torsion free $W(k)$-algebra $R_{\nu}$,
a representation $\rho_{\nu}: W_F \rightarrow \GL_n(R_{\nu})$, and an algebraic group $G_{\nu}$ with an action on $\Spec R_{\nu}$
that are ``universal'' (in a sense made precise in~\cite{curtis}, Proposition 9.2) for suitably rigidified representations of $W_F$ whose restriction to
$I_F^{(\ell)}$ is isomorphic to $\nu$.

In section 10 of~\cite{curtis} we formulate two precise conjectures.  First, we have the following conjecture, which we will call
the ``weak conjecture'':

\begin{conj}[\cite{curtis}, Conjecture 10.2] \label{conj:weak}
There is a map $\LL_{\nu}: eZ_n \rightarrow R_{\nu}$ ``compatible with local Langlands'' in the sense that for any $x: R_{\nu} \rightarrow \overline{\CK}$,
the composition with $\LL_{\nu}$ is the map $eZ_n \rightarrow \overline{\CK}$ giving the action of $eZ_n$ on the representation
$\Pi_x$ corresponding to the specialization $\rho_{\nu,x}$ of $\rho_{\nu}$ at $x$ via local Langlands.
\end{conj}

The map $\LL_{\nu}$ is unique if it exists.  Moreover, because the representation $\rho_{\nu}$ is $G_{\nu}$-invariant, the
image of $\LL_{\nu}$ is contained in the subring $R_{\nu}^{\inv}$ of $G_{\nu}$-invariant elements of $R_{\nu}$.  We then
make the further ``strong'' conjecture:

\begin{conj}[\cite{curtis}, Conjecture 10.3] \label{conj:strong}
The map $\LL_{\nu}: eZ_n \rightarrow R_{\nu}$ of the ``weak conjecture'' identifies $eZ_n$ with $R_{\nu}^{\inv}$.
\end{conj}

The weak conjecture is easy to verify after inverting $\ell$.  That is, Theorem 10.4 of~\cite{curtis} constructs a map:
$$\LL_{\nu}[\frac{1}{\ell}]: eZ_n \rightarrow R_{\nu}[\frac{1}{\ell}]$$
compatible with local Langlands; the weak conjecture amounts to showing that $\LL_{\nu}[\frac{1}{\ell}](eZ_n)$ is contained
in $R_{\nu}$.  In fact, by Theorem 10.5 of~\cite{curtis}, the image of $eZ_n$ under $\LL_{\nu}[\frac{1}{\ell}]$ of $eZ_n$
is contained in the normalization $\tR_{\nu}$ of $R_{\nu}$.  We denote by $\tLL_{\nu}$ the resulting map $eZ_n \rightarrow \tR_{\nu}$.

The goal of this section is to prove both of these conjectures for all $n$.  The argument will be inductive, and makes use
of the following key result, which is the main theorem of~\cite{curtis}:

\begin{thm}[\cite{curtis}, Theorem 11.1]
Suppose that the ``weak conjecture'' holds for all $G_m$ with $m \leq n$.  Then the ``strong conjecture'' holds for all such $G_m$.
\end{thm}

Both the ``weak conjecture'' and the ``strong conjecture'' hold for $G_1$ (this is an easy consequence of local class field theory).
Thus to prove both the weak and the strong conjecture for $G_n$, it suffices to show:

\begin{thm} \label{thm:induct}
Suppose the ``strong conjecture'' holds for $G_{n-1}$.  Then the ``weak conjecture'' holds for $G_n$.
\end{thm}

The main ingredients in the proof of this theorem are Corollary~\ref{cor:main descent},
together with a theory of $\gamma$-factors for families of
representations of $W_F$.  We first recall the latter theory, which was developed in~\cite{galoisgamma}.

If $\kappa$ is a field of characteristic zero containing $W(k)$, and $\rho: W_F \rightarrow \GL_n(\kappa)$ is a Weil group representation,
then there is a rational function $\gamma(\rho, X, \psi)$ in $\kappa(X)$, called the Deligne--Langlands $\gamma$-factor
of $\rho$.  This $\gamma$-factor is compatible with the local Langlands correspondence, in the sense that if
$\pi$ is an absolutely irreducible representation of $G_n$ over $\kappa$, and $\pi'$ is an absolutely irreducible representation
of $G_{n-1}$ over $\kappa$, then the $\gamma$-factor $\gamma(\pi \times \pi', X, \psi)$ of the pair $(\pi,\pi')$ coincides with
the $\gamma$-factor $\gamma(\rho \otimes \rho', X, \psi)$, where $\rho$ and $\rho'$ correspond to $\pi$ and $\pi'$
via local Langlands.

The main result of~\cite{galoisgamma} extends this construction to families of Weil representations.  In particular, one has:
\begin{thm}[\cite{galoisgamma}, Theorem 1.1] \label{thm:local galois gamma}
Let $R$ be a Noetherian $W(k)$-algebra and let $\rho: W_F \rightarrow \GL_n(R)$ be a representation that is
$\ell$-adically continuous in the sense of~\cite{galoisgamma}, section 4.
Then there exists an element $\gamma_R(\rho, X, \psi)$ of $S^{-1}R[X,X^{-1}]$ with the following properties:
\begin{enumerate}
\item If $f: R \rightarrow R'$ is a local homomorphism of complete Noetherian local $W(k)$-algebras, then one has:
$$f(\gamma_R(\rho, X, \psi)) = \gamma_{R'}(\rho \otimes_R R', X, \psi),$$
where we have extended $f$ to a map $S^{-1}R[X,X^{-1}] \rightarrow (S')^{-1}R'[X,X^{-1}]$ in the obvious way.
\item If $R$ is a field of characteristic zero, then $\gamma_R(\rho, X, \psi)$ coincides with the Deligne-Langlands
gamma factor $\gamma(\rho,X,\psi).$
\end{enumerate}
\end{thm}

Note that if $R$ is reduced and $\ell$-torsion free then the second property characterizes $\gamma_R(\rho,X,\psi)$ uniquely; in particular
it pins down the ``universal $\gamma$-factors'' $\gamma(\rho_{\nu},X,\psi)$ for all $\nu$.  Since any $\rho$ arises by base change from
some finite collection of these universal $\gamma$-factors, the two properties uniquely characterize the association
$\rho \mapsto \gamma(\rho,X,\psi)$.

{\em Proof of Theorem~\ref{thm:induct}:} Suppose that Conjecture~\ref{conj:strong} holds for $G_{n-1}$,
and fix a primitive idempotent $e$ of $Z_n$, corresponding to a representation $\nu$ of $I_F^{(\ell)}$.
We have a map: $\tLL_{\nu}: eZ_n \rightarrow \tR_{\nu}$.  Let $Z'$ be the preimage of $R_{\nu}$
under the map $\tLL_{\nu}$; it suffices to show that $Z' = eZ_n$.  Equivalently, it suffices to show
that the map giving the action of $eZ_n$ on $eW_n$ factors through $Z'$ (as this map is the identity on
$eZ_n$.)  

On the other hand, since $\tR_{\nu}$ is finitely generated as an $R_{\nu}$-module, $eZ_n$ is finitely generated
as a $Z'$-module.  Thus, by Corollary~\ref{cor:main descent}, it suffices to show, for each primitive idempotent
$e'$ of $Z_{n-1}$, that $\gamma(eW_n \times e'W_{n-1}, X^{-1}, \psi)$ and $\gamma(eW_n^{\iota} \times e'W_{n-1}, X, \psi)$
have coefficients in $Z' \otimes e'Z_{n-1}$.

Let $\nu'$ be the representation of $I_F^{(\ell)}$ corresponding to $e'$, and consider the map:
$$\tLL_{\nu} \otimes \LL_{\nu'}: eZ_n \otimes e'Z_{n-1} \rightarrow \tR_{\nu} \otimes R_{\nu'}^{\inv}.$$
By Conjecture~\ref{conj:strong} for $G_{n-1}$, the map $\LL_{\nu'}$ is an isomorphism of $e'Z_{n-1}$
with $R_{\nu'}^{\inv}$, so $Z' \otimes e'Z_{n-1}$ is the preimage of $R_{\nu} \otimes R_{\nu'}^{\inv}$
under $\tLL_{\nu} \otimes \LL_{\nu'}$.

Since this map is compatible with the local Langlands correspondence, the image of $\gamma(eW_n \times e'W_{n-1}, X^{-1}, \psi)$
under this map is the gamma factor $\gamma_{R_{\nu} \otimes{R_{\nu'}}}(\rho_{\nu} \otimes \rho_{\nu'}, X^{-1},\psi)$.
In particular it has coefficients in the tensor product
$R_{\nu} \otimes R_{\nu'}^{\inv}$, and so $\gamma(eW_n \times e'W_{n-1}, X^{-1},\psi)$
has coefficients in $Z' \otimes e'Z_{n-1}$ as claimed.  The proof that $\gamma(eW_n^{\iota} \times e'W_{n-1}, X,\psi)$
has coefficients in $Z' \otimes e'Z_{n-1}$ is similar.  $\Box$

The inductive argument given above then gives:

\begin{cor} Both Conjecture~\ref{conj:weak} and Conjecture~\ref{conj:strong} hold for all $n$.
\end{cor}

As an immediate consequence, we deduce that the elements of the Bernstein center constructed in~\cite{chenevier}, Proposition 3.11,
are integral:

\begin{cor} Let $w$ be an element of $W_F$.  Then there exists an element $z_w$ of $Z_n$ such that, for all irreducible smooth
$\overline{\QQ}_{\ell}$-representations $\pi$ of $\GL_n$, the action of $z_w$ on $\pi$ is given by $\tr \rho(w)$, where 
$\rho: W_F \rightarrow \GL_n(\overline{\QQ}_{\ell})$ corresponds to $\pi$ via local Langlands.
\end{cor}
\begin{proof}
For each $\nu$, consider $\tr \rho_{\nu}(w)$.  This is an element of $R_{\nu}^{\inv}$, and hence corresponds to an element of
$e_{\nu} Z_n$, which we denote by $z_{w,\nu}$.  The desired element $z_w$ is the product, over all $\nu$, of the $z_{w,\nu}$.
\end{proof}

\end{document}